\documentclass[12pt]{amsart}

\usepackage{amssymb, amsthm, amsmath, graphicx, yfonts, wasysym, appendix, url, verbatim}
\usepackage[margin=3cm]{geometry}
\usepackage{amsbsy}

\usepackage{amssymb}
\usepackage{amscd}
\usepackage{hyperref}
\usepackage{mathtools}
\usepackage{mathrsfs}

\def\Cal{\mathcal}

\let\phi\varphi

\newcommand{\defn}[1]{\emph{#1}}
\DeclareMathOperator{\dvol}{dV}

\DeclareMathOperator{\Conf}{Conf}
\DeclareMathOperator{\Ric}{Ric}
\DeclareMathOperator{\tr}{tr}
\newcommand{\suchthat}{\mathrel{}:\mathrel{}}

\newcommand{\mE}{\mathcal{E}}
\newcommand{\mG}{\mathcal{G}}
\newcommand{\cmE}{\widetilde{\mathcal{E}}}
\newcommand{\cmG}{\widetilde{\mathcal{G}}}
\newcommand{\cdelta}{\widetilde{\delta}}
\newcommand{\cDelta}{\widetilde{\Delta}}
\newcommand{\cnabla}{\widetilde{\nabla}}
\newcommand{\cg}{\widetilde{g}}
\newcommand{\cu}{\widetilde{u}}
\newcommand{\hg}{\widehat{g}}

\newcommand{\cc}{\mathfrak{c}}

\newcommand{\bg}{\boldsymbol{g}}

\newcommand{\ce}{\mathcal{E}}

\newcommand{\rpl}                         % +) or <+
{\mbox{$
\begin{picture}(12.7,8)(-.5,-1)
\put(0,0.2){$+$}
\put(4.2,2.8){\oval(8,8)[r]}
\end{picture}$}}

\numberwithin{equation}{section}

\newtheorem{thm}{Theorem}[section]
\newtheorem*{thm*}{Theorem \thesubsection}

\newtheorem*{lemma*}{Lemma \thesubsection}
\newtheorem*{prop*}{Proposition \thesubsection}
\newtheorem*{cor*}{Corollary \thesubsection}
\newtheorem{conj}[thm]{Conjecture}

\theoremstyle{definition}

\newtheorem*{definition*}{Definition \thesubsection}

\newtheorem*{example*}{Example \thesubsection}
\theoremstyle{remark}

\newtheorem*{remark*}{Remark \thesubsection}

\def\sideremark#1{\ifvmode\leavevmode\fi\vadjust{\vbox to0pt{\vss% the remark
 \hbox to 0pt{\hskip\hsize\hskip1em%                          will appear only
 \vbox{\hsize3cm\tiny\raggedright\pretolerance10000%          on the side
  \noindent #1\hfill}\hss}\vbox to8pt{\vfil}\vss}}}%
\begin{document}
\renewcommand{\today}{}
\title{The GJMS operators in geometry, analysis, and physics}

\author{Jeffrey S.\ Case}
\address{J.S.C:Department of Mathematics \\ Penn State University \\ University Park, PA 16802, USA}
\email{jscase@psu.edu}
\author{A.\ Rod Gover}
\address{A.R.G.:Department of Mathematics\\
  The University of Auckland\\
  Private Bag 92019\\
  Auckland 1142\\
  New Zealand}
\email{r.gover@auckland.ac.nz}

\begin{abstract}
 The GJMS operators, introduced by Graham, Jenne, Mason, and Sparling, are a family of conformally invariant linear differential operators with leading term a power of the Laplacian.
 These operators and their method of construction have had a major impact in geometry, analysis, and physics.
 We describe the GJMS operators and their construction, and briefly survey their importance and impact.
\end{abstract}

% This is because the editor on my iPad doesn't recognize the optional argument
% [2020] to \subjclass
\makeatletter
\@namedef{subjclassname@2020}{\textup{2020} Mathematics Subject Classification}
\makeatother

\maketitle

\pagestyle{myheadings} \markboth{Case \& Gover}{The GJMS operators in geometry, analysis, and physics}

\section{Introduction}

Riemannian geometry is the study of (pseudo-)Riemannian manifolds up to isometry.
Natural linear differential operators on Riemannian manifolds, which by definition commute with isometries, play a fundamental role in mathematics and physics.
One of the most important such operators is the Laplacian $\Delta$;
see~\cite{Lablee2015}.

Conformal geometry is the study of Riemannian manifolds up to conformal equivalence, where two metrics $g$ and $\hg$ are {\em conformally equivalent} if $\hg = e^{2\Upsilon}g $ for some smooth function $\Upsilon$.
On surfaces the Laplacian is conformally invariant;
more precisely,
\begin{equation}
 \label{eqn:surface-transformation}
 e^{2\Upsilon}\Delta^{\hg} = \Delta^g \, .
\end{equation}
The GJMS operators of Graham, Jenne, Mason, and Sparling, published in the landmark paper~\cite{GJMS} in JLMS in 1992, gave an essentially complete answer to the question of the existence of natural, conformally invariant, linear differential operators with leading-order term a power of the Laplacian.
The goals of this survey are to describe the GJMS construction and to highlight a variety of its impacts in mathematics and physics.
We begin with a brief account of the historical background.

It was realised by Bateman~\cite{Bateman1909} and
Cunningham~\cite{Cunningham1910} in 1909 that the (indefinite) Laplacian, or wave operator, and Maxwell's equations, respectively, are conformally invariant in Minkowski four-space.  More
generally for the Laplacian, if $\Phi$ is a diffeomorphism of (an open subset of)
$\mathbb{R}^n$ such that $\Phi^\ast g = e^{2\Upsilon}g$, where $g$ is
the standard flat metric, then
\begin{equation}
 \label{eqn:flat-conformal-laplacian}
 \Delta^{\Phi^\ast g}u = e^{-\frac{n+2}{2}\Upsilon} \Delta^g \bigl(e^{\frac{n-2}{2}\Upsilon}u\bigr)
\end{equation}
for all smooth functions $u$ in (the same open subset of) $\mathbb{R}^n$.
The conformal transformation law~\eqref{eqn:flat-conformal-laplacian} is not valid on general manifolds.
However, the \defn{conformal Laplacian}
\begin{equation}
 \label{P2expl}
 P_2 := \Delta + \frac{n-2}{4(n-1)}R ,
\end{equation}
obtained by adding a multiple of the scalar curvature $R$ as a zeroth-order operator to the Laplacian (with sign convention $\Delta \geq 0$ in positive signature), is conformally invariant on $n$-manifolds:
if $\hg = e^{2\Upsilon}g$, then
\begin{equation}
 \label{P2trans}
 P_2^{\hg} = e^{-\frac{n+2}{2}\Upsilon} \circ P_2^g \circ e^{\frac{n-2}{2}\Upsilon} ,
\end{equation}
where the exponentials are regarded as multiplication operators.
This generalises Equations~\eqref{eqn:surface-transformation} and~\eqref{eqn:flat-conformal-laplacian}.
Equivalently, $P_2$ defines a conformally invariant operator
\begin{equation*}
 P_2 \colon \ce\left( -\frac{n-2}{2} \right) \to \ce\left( -\frac{n+2}{2} \right) ,
\end{equation*}
where $\ce(w)$ is the set of \defn{conformal densities} of weight $w$;
i.e.\ equivalence classes of $(g,u) \sim (e^{2\Upsilon}g,e^{w\Upsilon}u)$ on pairs of a metric and a function.
The conformal transformation law for the scalar curvature was known to Schouten~\cite{Schouten1921} in 1921, but its expression using a \emph{linear} operator was popularised by Yamabe~\cite{Yamabe1960} in 1960.

The existence of the conformal Laplacian raises the question of
whether there are analogous higher-order operators.
In 1977, Jakobsen and Vergne~\cite{JakobsenVergne1977} made a major step in this direction by
proving that integer powers of the Laplacian are conformally invariant, in a
sense analogous to Equation~\eqref{eqn:flat-conformal-laplacian}, on
Minkowksi four-space; this generalises to all dimensions~\cite{GJMS}.
This led many authors to
seek, in analogy with Equation~\eqref{P2expl}, lower-order correction
terms that yield a conformally invariant operator with leading-order
term $\Delta^k$ on general manifolds.  When $k=2$, this goal was
achieved by several
authors~\cite{EastwoodSinger1985,FradkinTseytlin1982,Paneitz1983,Riegert1984},
including in all dimensions $n \not= 1,2$ by Paneitz in 1983, and the
resulting $P_4$ is now called the {\em Paneitz operator}.  When $k=3$, this
goal was achieved in all dimensions $n \not= 1,2,4$ by
W\"unsch~\cite{Wunsch1986} in 1985.  The approach of these
papers was to solve an undetermined coefficients problem whose
computational complexity quickly gets out of hand.

Graham, Jenne, Mason, and Sparling~\cite{GJMS} introduced a new conceptual insight which allowed them to completely bypass this computational complexity and construct higher-order conformally invariant operators with leading-order term $\Delta^k$:

\begin{thm}
 \label{GJMSthm}
 Fix an integer $n \geq 3$.
 Let $k \in \mathbb{N}$;
 if $n$ is even, then assume additionally that $k \leq n/2$.
 There is a natural, conformally invariant, linear differential operator
 \begin{equation*}
  P_{2k} \colon \ce\left( -\frac{n-2k}{2} \right) \to \ce\left( -\frac{n+2k}{2} \right)
 \end{equation*}
 on $n$-manifolds with leading-order term $\Delta^k$.
\end{thm}

Equivalently, Theorem~\ref{GJMSthm} constructs natural linear differential operators $P_{2k}^g$ with leading-order term $\Delta^k$ such that if $\hg = e^{2\Upsilon}g$, then
\begin{equation}
 \label{P2ktrans}
 P_{2k}^{\hg} = e^{-\frac{n+2k}{2}\Upsilon} \circ P_{2k}^g \circ e^{\frac{n-2k}{2}\Upsilon} .
\end{equation}

Remarkably, Theorem~\ref{GJMSthm} is sharp: Graham~\cite{Graham1992}
showed, in the same issue of JLMS as~\cite{GJMS}, that the operator $P_6$
cannot exist on general conformal $4$-manifolds, and Gover and
Hirachi~\cite{GoverHirachi2004} proved that the operators $P_{2k}$
cannot exist when $k > n/2$ and $n$ is even.

The significant new conceptual insight in~\cite{GJMS} is that the operators $P_{2k}$ have
a simple \emph{holographic} construction; i.e.\ they are
determined by natural operators in a canonically associated Riemannian
manifold.  In this case, $P_{2k}$ is determined by the Laplace powers
$\widetilde{\Delta}^k$ in the Fefferman--Graham ambient
space~\cite{FeffermanGraham1985,FeffermanGraham2012}.
This idea was used by Bateman~\cite{Bateman1909} to prove the
conformal invariance of the Laplacian in Minkowski four-space (see
also~\cite{Dirac1936,HughstonHurd1983}), with the conformal Laplacian itself
recovered by Fefferman and Graham~\cite{FeffermanGraham1985} in 1985.
In fact, there are two different, but equivalent, constructions of the
GJMS operators, which 
are powerful enough to allow one to
prove additional key properties of the GJMS operators: they are
formally self-adjoint~\cite{GrZ}, they give rise to scalar invariants that generalise the
Gauss and scalar curvatures \cite{Branson1995}, and they factor as compositions of
second-order Schr\"odinger operators at Einstein
manifolds~\cite{FeffermanGraham2012,G-06}.  These constructions have
been adapted to produce conformally invariant operators in a variety
of settings, including on sections of other tensor
bundles~\cite{BransonGover2005,Matsumoto2013}, on spin bundles~\cite{EelbodeSoucek2010,Fischmann2014,FischmannKrattenthalerSomberg2015,HollandSparling2001}, on
submanifolds~\cite{CGK25,GoverWaldron2021}, and on CR
manifolds~\cite{GoGr,Hirachi2013}; they have also been generalised to
produce nonlocal linear conformally invariant pseudodifferential
operators with leading-order term a fractional power of the
Laplacian~\cite{BransonGover2001,ChangGonzalez2011,GrZ}.

The GJMS operators and their holographic method of construction continue to inspire research in many directions.
We indicate some of these by highlighting a few open conjectures.
This survey is organised as follows:

In Section~\ref{sec:gjms-construction} we explain the two holographic constructions of the GJMS operators.
We then describe a related construction on Poincar\'e spaces via scattering theory~\cite{GrZ}, which has been highly influential in the context of the AdS/CFT correspondence~\cite{GubserKlebanovPolyakov1998,Maldacena1998,Witten1998} in physics.
We conclude by discussing local formulae for the GJMS operators. 

In Section~\ref{sec:alg} we discuss the aforementioned properties of the GJMS operators and their applications.
We also discuss related constructions of invariant operators.

In Section~\ref{sec:ga} we discuss applications of the GJMS operators
to problems in geometric analysis related to the rigidity and
stability of Sobolev inequalities and to various analogues of the
Yamabe Problem.

\section{The ambient space and the GJMS construction} \label{sec:gjms-construction}

Let $(M,\cc)$ be a conformal $n$-manifold.
This data is equivalent to a ray subbundle
\begin{equation*}
 \mathcal{G} := \left\{ (x,g_x) \suchthat x \in M, g \in \mathfrak{c} \right\} \subset S^2T^\ast M .
\end{equation*}
The projection $\pi(x,g_x) := x$ and dilation $\delta_s(x,g_x) := (x,s^2g_x)$ make $\pi \colon \mG \to M$ into a principal $\mathbb{R}_+$-bundle.
There is a tautological symmetric $(0,2)$-tensor $\bg$ defined by
\begin{equation*}
 \bg(X,Y) := g_x(\pi_\ast X,\pi_\ast Y)
\end{equation*}
for $X,Y \in T_{(x,g_x)}\mathcal{G}$.
The set of \defn{conformal densities} of weight $w \in \mathbb{R}$ is
\begin{equation*}
 \mE(w) := \left\{ u \in C^\infty(\mG) \suchthat \delta_s^\ast u = s^w u \right\} .
\end{equation*}
Evidently $u \in \mE(w)$ if and only if $u^g(x) := u(x,g_x)$ satisfies
\begin{equation*}
 (g,u^g) \sim (e^{2\Upsilon}g,e^{w\Upsilon}u^{e^{2\Upsilon}g}) ,
\end{equation*}
so that this definition of $\mE(w)$ agrees with its definition in the introduction.

The (Fefferman--Graham) \defn{ambient space} $(\cmG,\cg)$ for
$(M,\cc)$ is a smooth Riemannian $(n+2)$-manifold
with a proper embedding $\iota \colon \mG \to \cmG$ and
dilations $\cdelta_s \colon \cmG \to \cmG$ such that
\begin{enumerate}
 \item $\iota^\ast\cg = \bg$,
 \item $\cdelta_s^\ast\cg = s^2\cg$,
 \item $\iota \circ \delta_s = \cdelta_s \circ \iota$ for all $s \in \mathbb{R}_+$, and
 \item $\Ric(\cg) = O^+(\sigma^{(n-2)/2})$, if $n \geq 4$ is even, and $\Ric(\cg) = O(\sigma^\infty)$ otherwise.
\end{enumerate}
Here $\sigma$ is a defining function for $\iota(\mG)$ which is
homogeneous of degree $2$---i.e.\ $\cdelta_s^\ast\sigma =
s^2\sigma$---and $O^+(\sigma^m)$ is the set of all symmetric
$(0,2)$-tensor fields $T$ in $O(\sigma^m)$ such that for each $z = (x,g_x) \in \mG$,
there is a $t \in S^2T_{x}^\ast M$ such that
$\iota^\ast(\sigma^{-m}T)_z = (\pi^\ast t)_z$ and $\tr_{g_x}t = 0$.
Fefferman and Graham~\cite{FeffermanGraham2012} showed that when $n
\geq 3$, ambient spaces exist and are unique up to shrinking $\cmG$,
applying a $\cdelta_s$-equivariant diffeomorphism which restricts to
the identity on $\iota (\mG)$, and modifying $\cg$ by adding a term in
$O(\sigma^\infty)$ or $O^+(\sigma^{n/2})$ if $n$ is odd or even,
respectively.
When $(M,\mathfrak{c})$ is the
$n$-sphere with its standard conformal structure, one can take
$(\cmG,\cg)$ as $(n+2)$-dimensional Minkowski space, with $\mG$ the null cone therein and dilation
arising from scalar multiplication.
This case provides the flat model that the ambient space generalises. 
When $n \geq 4$ is even, the \defn{Fefferman--Graham obstruction
  tensor}
\begin{equation*}
 \mathcal{O} := c_n\iota^\ast(\sigma^{-(n-2)/2}\Ric(\cg))
\end{equation*}
determines a conformally invariant, trace-free, divergence-free
symmetric $(0,2)$-tensor that obstructs the existence of an ambient
metric $\cg$ satisfying $\Ric(\cg) = O(\sigma^\infty)$.
In the case of
dimension $n=4$, the obstruction tensor $ \mathcal{O}$ is the Bach tensor~\cite{Bach1921}.

The ambient space was initially used to construct scalar conformal invariants.
Set
\begin{equation*}
	\cmE(w) := \bigl\{ \cu \in C^\infty(\cmG) \suchthat \cdelta_s^\ast\cu = s^w\cu \bigr\}
\end{equation*}
and note that $\iota^\ast$ maps $\cmE(w)$ onto $\mE(w)$.
It follows that if $\widetilde{I} \in \cmE(w)$ is a natural scalar Riemannian invariant of the ambient space, then $\iota^\ast\widetilde{I}$ defines a scalar conformal invariant provided it is independent of the ambiguity of the ambient metric.
See~\cite{BaileyEastwoodGraham1994,Gover2001} for applications to the (partial) classification of scalar conformal invariants.

The key insight of Graham, Jenne, Mason, and Sparling~\cite{GJMS} is
that powers of the Laplacian $\cDelta$ of $\cg$ determine conformally
invariant operators on density
bundles. These are the \emph{GJMS operators}.  They in fact gave two equivalent constructions of these
operators, both of which are very important.

Their first construction establishes the surprising result that
$\cDelta^k \colon \cmE(w) \to \cmE(w-2k)$ descends to $\mG$---in the
sense that $\iota^\ast(\cDelta^k\cu)$ depends only on
$\iota^\ast\cu$---if and only if $w = -\frac{n-2k}{2}$.  Thus one can
define
\begin{equation}
 \label{eqn:gjms-defn}
 P_{2k}u := \iota^\ast\bigl( \cDelta^k\cu \bigr) \in \mE\left( -\frac{n+2k}{2} \right)
\end{equation}
for some, and hence any, homogeneous extension $\cu$ of $u \in
\mE\bigl(-\frac{n-2k}{2}\bigr)$.
Their proof uses a certain defining function $\sigma \in \cmE(2)$ for $\iota(\mG) \subset \cmG$ with the property that that, as operators, $\sigma$, $\cDelta$, and $[\cDelta,\sigma]$ generate the Lie algebra $\mathfrak{sl}(2)$.
The result then follows from a standard $\mathfrak{sl}(2)$ identity.
The restriction $k \leq n/2$ in
Theorem~\ref{GJMSthm} when $n$ is even guarantees that $P_{2k}$ is
independent of the $O^+(\sigma^{n/2})$-ambiguity of the ambient metric
$\cg$, and hence is well-defined.

Their second construction exhibits the GJMS operators as the obstruction to finding a smooth, harmonic, homogeneous extension $\cu$ of $u \in \mE\bigl(-\frac{n-2k}{2}\bigr)$.
More precisely, they show that there is an extension $\cu$, unique modulo $O(\sigma^k)$, such that $\cDelta\cu = O(\sigma^{k-1})$, and that one can take $\cDelta\cu = O(\sigma^\infty)$ if and only if the obstruction
\begin{equation*}
	P_{2k}^{(\mathrm{obs})}u := \iota^\ast\bigl( \sigma^{1-k}\cDelta\cu \bigr)
\end{equation*}
vanishes.
Moreover, they show that $P_{2k}^{(\mathrm{obs})}=c_kP_{2k}$ for an explicit nonzero constant $c_k$.

The second construction of the GJMS operators has a useful interpretation in terms of scattering theory.
Consider the restriction $g_+$ of $\cg$ to the set $X := \{ \sigma = -1 \}$.
Then $(X,g_+)$ is a Poincar\'e manifold~\cite{FeffermanGraham2012};
i.e.\ it is asymptotically hyperbolic and formally solves $\Ric(g_+) + ng_+ = 0$ to a dimension-dependent order.
One directly computes that if $\cu \in \cmE\bigl(-\frac{n-2k}{2}\bigr)$, then $\cDelta\cu = 0$ if and only if
\begin{equation}
	\label{eqn:scattering}
	\Delta_{g_+}U - s(n-s)U = 0 ,
\end{equation}
where $U := \cu \rvert_X$ and $s := \frac{n}{2} + k$.
Graham and Zworski~\cite{GrZ} proved that for generic $s>n/2$, if $r$ is a defining function for a suitably defined boundary $\partial X$ of $X$ and if $f \in C^\infty(\partial X)$, then there is a unique solution $U$ of~\eqref{eqn:scattering} and, moreover,
\begin{equation*}
 U = Fr^{n-s} + Gr^s 
\end{equation*}
for $F,G$ sufficiently smooth functions on $\overline{X}$ with $F\rvert_{r=0}=f$.
The interpretation of GJMS operators as obstructions implies that the \defn{scattering operator}
\begin{equation}
 \label{eqn:scattering-matrix}
 S(s)f := G\rvert_{r=0}
\end{equation}
has a pole at each $s = \frac{n}{2}+k$, $k \in \mathbb{N}$, and this pole is proportional to $P_{2k}f$.

One difficulty with the definition~\eqref{eqn:gjms-defn} is that producing a local formula for $P_{2k}$ requires computing the asymptotic expansion of the ambient metric $\cg$.
The nonlinearity of the equation $\Ric(\cg)=0$ makes this quite challenging.
Fefferman and Graham~\cite{FeffermanGraham2013} used the first construction to give a direct proof of the holographic formula
\begin{equation}
 \label{eqn:juhl-formula}
 P_{2k} = \sum_{\lvert I \rvert = k} n_I \mathcal{M}_{I} ,
\end{equation}
originally due to Juhl~\cite{Juhl2013}.
Here $\mathcal{M}_{I} = \mathcal{M}_{i_1} \dotsm \mathcal{M}_{i_k}$ is a composition of second-order operators appearing in the asymptotic expansion of $\cDelta$.
This decomposition is useful when studying Yamabe-type constants;
see Section~\ref{sec:GA/yamabe}.

Another route to producing local formulae for the GJMS operators is
through the tractor calculus~\cite{BaileyEastwoodGover1994}.
Although this
calculus is easily defined directly,
it can be recovered via the ambient metric.  The
\defn{standard tractor bundle} $\mathcal{T}[w]$ of weight $w$ can be
identified with the set of equivalence classes of vectors in the
pullback $\iota^\ast T\cmG$ of the tangent bundle $T\cmG$
that are homogeneous of degree
$w-1$ along the orbits of dilation~\cite{CapGover2003}.
Roughly
speaking, a tractor bundle $\mathcal{T}^\Phi[w]$ is a bundle obtained
by taking tensor powers of standard tractor bundles and their
duals. Then the standard \defn{tractor connection} arises from the
pullback of $\cnabla$ to $\mG$, giving a conformally invariant
connection on $\mathcal{T}^\Phi[0]$.  Certain linear combinations of
$\widetilde{\nabla}$ and $\widetilde{T}\cDelta$, where $\widetilde{T}$ is the infinitesimal generator of dilations,
produce conformally invariant maps $\mathcal{T}^\Phi[w] \to \mathcal{T} \otimes
\mathcal{T}^\Phi[w-1]$ called \defn{tractor-$D$ operators}.
Gover and Peterson~\cite{GoverPeterson2003} used this to produce formulae
for the GJMS operators that involve only conformally invariant
operators between tractor bundles;
these formulae do not require the asymptotic expansion of the ambient metric.
While these formulae are also quite complicated in general, they are
helpful in deriving some of the properties and generalisations of
Section~\ref{sec:alg}.

\section{Properties and generalisations} \label{sec:alg}

  \subsection{$Q$-curvatures}
  \label{sec:alg/q}

One particularly nice application of the GJMS operators is to the
construction and study of the $Q$-curvatures.  These invariants, first
introduced by Branson and {\O}rsted~\cite{BransonOrsted1991} in
dimension four and by Branson~\cite{Branson1995} in general dimension,
arise out of the fundamental properties of the GJMS operators, as can
already been seen by considering the conformal Laplacian $P_2$.  On
the one hand, the explicit formula~\eqref{P2expl} implies that $P_2$
is formally self-adjoint and that $P_2(1)=0$ in dimension two.
On the other hand, the conformal transformation
formula~\eqref{P2trans} can be rewritten as
\begin{equation}\label{dimc}
 \frac{2}{n-2}e^{\frac{n+2}{2}\Upsilon}P_2^{e^{2\Upsilon}g}(1) = \frac{2}{n-2}P_2^g(1) + \frac{2}{n-2}P_2^g(e^{\frac{n-2}{2}\Upsilon}-1) 
\end{equation}
in dimension $n \geq 3$.
Since $\frac{2}{n-2}P_2^g(1) = \frac{1}{2(n-1)}R^g$, formally taking the limit $n \to 2$ recovers the transformation law
\begin{equation}
 \label{preQ}
 e^{2\Upsilon}K^{e^{2\Upsilon}g} = K^g + \Delta^g\Upsilon
\end{equation}
in dimension two, where $K^g = \frac{1}{2}R^g$ is the Gauss curvature.
Integrating Equation~\eqref{preQ} implies that the total Gauss curvature is a conformal invariant of compact surfaces.
Of course, $\int K^g \dvol_g = 2\pi\chi(M)$ in dimension two, so it is in fact a topological invariant.

Since constants are in $\widetilde{\mathcal{E}}(0)$, Equation~\eqref{eqn:gjms-defn} implies that $P_n(1)=0$ in even dimension~$n$.
Branson~\cite{Branson1995} further observed that the GJMS construction is such that one can formally define the \defn{$Q$-curvature} of a $2k$-manifold as the limit
\begin{equation*}
 Q^g := \lim_{n \to 2k} \frac{1}{n-2k}P_{2k}^g(1) ;
\end{equation*}
see~\cite{FeffermanGraham2002,FeffermanHirachi2003,GoverPeterson2003} for alternate constructions. Arguing as in the recovery of \eqref{preQ} from \eqref{dimc},
it then follows from Equation~\eqref{P2ktrans} that
\begin{equation}
 \label{Qtrans}
 e^{2k\Upsilon}Q^{e^{2\Upsilon}g} = Q^g + P_{2k}^g(\Upsilon) .
\end{equation}
Assuming that $P_{2k}$ is formally self-adjoint, one concludes that the total $Q$-curvature is a conformal invariant of compact $2k$-manifolds.
Formal self-adjointness was later proved by Graham and Zworski~\cite{GrZ}.

The total $Q$-curvature is not a topological invariant in general.
In dimension four,
\begin{equation*}
 8\pi^2\chi(M) = \int_M Q \dvol + \int_M \lvert W \rvert^2 \dvol ,
\end{equation*}
where $W$ is the Weyl tensor~\cite{BransonOrsted1991}.
More generally, Alexakis~\cite{Alexakis2012} showed that if $I$ is a scalar Riemannian invariant whose integral is conformally invariant on compact $n$-manifolds, then there is a constant $c_I$ and and a scalar conformal invariant $\mathcal{W}_I$ such that
\begin{equation}
 \label{eqn:alexakis}
 c_I\chi(M) = \int_M I \dvol + \int_M \mathcal{W}_I \dvol
\end{equation}
for any compact Riemannian $n$-manifold;
see~\cite{BoulangerFrancoisLazzarini2019} for an outline of a different proof.
Note that the scalar conformal invariant $\mathcal{W}_I$ is not uniquely determined in general~\cite{CKLTY25}.

\subsection{Relations to physics}
\label{subsec:physics}

The Lorentzian signature analogue of hyperbolic space is known as Anti-de Sitter (AdS) space.
Holographic constructions of conformal invariants via Poincar\'e manifolds provide an important motivation for the AdS/CFT correspondence~\cite{GubserKlebanovPolyakov1998,Maldacena1998,Witten1998};
i.e.\ the expected correspondence between a conformal field theory (CFT) on the boundary and the bulk geometry of AdS space.

On a Poincar\'e--Einstein manifold $(X^{n+1},g_+)$, the Einstein--Hilbert action $\int R \dvol$, which is proportional to the volume, is necessarily infinite.
For any geodesic defining function $r$ for $\partial X$---which uniquely determines a metric $g := r^2g_+\rvert_{T\partial X}$ on the conformal infinity---one computes that
\begin{equation}
 \label{eqn:volume-expansion}
 \int_{r > \varepsilon} \dvol_{g_+} = a_{(0)}\varepsilon^{-n} + a_{(2)}\varepsilon^{2-n} + \dotsm + \mathcal{L} \log \varepsilon^{-1} + \mathcal{V} + o(1)
\end{equation}
as $\varepsilon\to0$, where only the ``even'' terms
$a_{(2k)}\varepsilon^{2k-n}$ appear before the logarithmic term.  When
$n$ is odd, $\mathcal{L}$ vanishes and $\mathcal{V}$ is conformally
invariant.  When $n$ is even, $\mathcal{L}$ is conformally invariant
but generally nonzero, and $\mathcal{V}$ is not conformally invariant.
Instead,
\begin{equation}
 \label{eqn:anomaly}
 \left. \frac{d}{dt} \right|_{t=0} \mathcal{V}(e^{2t\Upsilon}g) = \int_M \Upsilon I_g \dvol_g 
\end{equation}
for some scalar invariant $I$; moreover, $\mathcal{L}(g) = \int I_g
\dvol_g$.  These results are due to Henningson and
Skenderis~\cite{HenningsonSkenderis1998}, who further observed that
$I$ calculates the conformal
anomaly~\cite{DeserDuffIsham1976,Hawking1977} of the the CFT on the
boundary.  It is known~\cite{GrZ} that $\mathcal{L}$ is proportional
to the total $Q$-curvature, leading to the result that the
Fefferman--Graham obstruction tensor $\mathcal{O}$ is proportional to
the metric gradient of $\int Q \dvol$; see~\cite{GrahamHirachi2005}.

In even dimensions, the logarithmic functional determinants of conformally invariant operators \emph{also} satisfy~\eqref{eqn:anomaly} for some invariant $I$~\cite{BransonOrsted1991dga,ParkerRosenberg1987}.
In dimensions $n \in \{ 2, 4, 6 \}$, formulae for the functional determinant are known~\cite{Branson1995,BransonOrsted1991,Polyakov1981};
they include the term
\begin{equation*}
 \mathcal{F}(u) = \int_{M^n} (uP_nu + 2Qu ) \dvol
\end{equation*}
(cf.\ Equation~\eqref{eqn:functional}).
The physics literature gives an alternative, holographic, route towards computing the functional determinants of GJMS operators using the realisation of the latter as poles of the scattering operator $S(s)$ and its realisation as the two-point function of a CFT function of conformal dimension $s$;
see~\cite{BuginiDiaz2019,Diaz2008} and references therein.

There are many names for scalar Riemannian invariants $I$ of homogeneity $-n$ with the property that there
is some functional $\mathcal{V}$ such that
Equation~\eqref{eqn:anomaly} holds;
e.g.\ \defn{conformal anomalies}~\cite{DeserSchwimmer1993},
\defn{conformal indices}~\cite{BransonOrsted1986}, or
\defn{conformally variational invariants}~\cite{CaseLinYuan2019}.  The
$Q$-curvature, Pfaffian, and scalar conformal invariants are
examples~\cite{Branson1995,DeserSchwimmer1993}.  A \defn{trivial}
anomaly is one obtained by taking $\mathcal{V}$ to be the integral of
a scalar Riemannian invariant of homogeneity $-n$; e.g.\ $\mathcal{V}
:= \frac{1}{2}\int J^2\dvol$ determines the trivial anomaly $\Delta J$
in dimension four.  Based on physical considerations, Deser and
Schwimmer~\cite{DeserSchwimmer1993} conjectured
(cf.\ \cite[Conjecture~14]{Branson2005}) that \emph{every} anomaly is
a linear combination of the Pfaffian, a scalar conformal invariant,
and a trivial anomaly.  Boulanger~\cite{Boulanger2007} sketched a
proof of this conjecture using BRST cohomology.

\subsection{Factorisation formulae} \label{factor}

Einstein metrics play a special role in conformal geometry, for
example through their distinguished ambient
metrics~\cite{FeffermanGraham2012} and their correspondence with
parallel tractors~\cite{GoverNurowski2006}.
This leads to a factorisation for the GJMS operators and a simple formula for the $Q$-curvature:

\begin{thm}
 \label{factorisation}
 Let $(M^n,g)$ be an Einstein manifold with $\Ric_g = (n-1)\lambda g$.
 Let $k$ be as in Theorem~\ref{GJMSthm}.
 Then
 \begin{equation}
  \label{eqn:factorisation}
  P_{2k} = \prod_{\ell=1}^k \left( \Delta^g + \frac{(n+2\ell-2)(n-2\ell)}{4}\lambda \right) .
 \end{equation}
 Moreover, if $n$ is even, then $Q = (n-1)!\lambda^n$.
\end{thm}

There are three different proofs of
Theorem~\ref{factorisation}.
One proof~\cite{G-06}
uses the correspondence between Einstein and parallel tractors to
express $P_{2k}$ as the composition $(I \cdot D)^k$, where $I$ is the
parallel tractor and $D$ is the tractor-$D$ operator.
Another proof~\cite{FeffermanGraham2012}
uses the special form of the Poincar\'e metric to explicitly solve Equation~\eqref{eqn:scattering}.
The third
proof, implicitly given in~\cite{Matsumoto2013}, uses a special choice of
extension to evaluate Equation~\eqref{eqn:gjms-defn}.

A first application of Theorem~\ref{factorisation} is to compute the
renormalised volume of an even-dimensional Poincar\'e--Einstein
manifold~\cite{ChangQingYang2006};
see~\cite{CGKTW25,GoverLatiniWaldronZhang25} for adaptations to other
settings.
The resulting formula is typically not explicit, in that it involves a
scalar conformal invariant that is determined indirectly by Alexakis'
decomposition theorem~\cite{Alexakis2012}.  An explicit formula has
recently been obtained using a generalisation of
Theorem~\ref{factorisation} to more general invariants of Einstein
manifolds~\cite{CKLTY25}.

A second application of Theorem~\ref{factorisation} is to scattering
theory.  For generic $\gamma \in (0,n/2)$, the rescaled scattering
operator $P_{2\gamma} = c_\gamma S\bigl(\frac{n}{2}+\gamma\bigr)$ of
Equation~\eqref{eqn:scattering-matrix} 
is a formally self-adjoint,
conformally covariant, pseudodifferential operator with leading-order
term $\Delta^\gamma$~\cite{GrZ}; these are the \defn{fractional GJMS
  operators}.  Using (a weighted analogue of)
Theorem~\ref{factorisation}, one can realise these operators as
Dirichlet-to-Neumann operators associated to the GJMS operators in the
interior, and also identify their Dirichlet
energies~\cite{CaseChang2016}.  In the special case of
$\mathbb{R}_+^{n+1}$ with the hyperbolic metric, this
recovers~\cite{ChangGonzalez2011} the energy approach to fractional
Laplacians popularised by Caffarelli and
Silvestre~\cite{CaffarelliSilvestre2007}.

A third application of Theorem~\ref{factorisation} is to a proof of the sharp Sobolev inequalities for the embeddings $W^{k,2} \hookrightarrow L^{2n/(n-2k)}$ on $\mathbb{R}^n$, $n > 2k$;
see Section~\ref{sec:GA/sobolev} for further discussion.
Crucially, this argument can be adapted to the CR GJMS operators of Section~\ref{CR-sect}.

Analogues of Theorem~\ref{factorisation} have also been proved for variants of the GJMS operators;
e.g.\ the GJMS operators on symmetric $2$-tensors~\cite{Matsumoto2013} and the GJMS operators on spinors~\cite{EelbodeSoucek2010,FischmannKrattenthalerSomberg2015}.
See below for additional examples.

\subsection{GJMS and related operators on differential forms}\label{forms}

\newcommand{\cC}{{\Cal C}}
\newcommand{\cd}{\widetilde{d}}
\renewcommand{\d}{\delta}

The GJMS constructions have inspired many
similar constructions in other
settings.  One particularly interesting construction, due to Branson and
Gover~\cite{BransonGover2005}, is to GJMS operators on differential
forms.  Since functions are $0$-forms, these are true generalisations
of the GJMS operators, and the behavior of the ``$Q$-curvature'' on
forms is particularly illuminating.

\begin{thm}
 \label{long} Fix a dimension $n$ and let $k \in \{ 0, 1, \dotsc,
 n \}$.  Let $\ell$ be a nonnegative integer; if $n$ is even, assume
 additionally that $\ell < n/2$, or that $k=0$ and $\ell =n/2$.  There is a natural, formally
 self-adjoint, conformally invariant, linear differential
 operator
 \begin{equation*}
  L_k^\ell \colon \ce^k(k+\ell-n/2) \to \ce^k (k-\ell-n/2)
 \end{equation*}
 of order $2\ell$ with
 leading-order term $(n-2k+2\ell)(\d d)^\ell +
 (n-2k-2\ell)(d\d)^\ell$, up to  a non-zero constant multiple.  Moreover, if $n$ is even and $k \leq n/2$,
 then
 \begin{equation*}
  L_k^{n/2-k} = (n-2k+2)\d Q_{k+1}d
 \end{equation*}
 for $Q_k$ a natural linear differential operator
 such that if $\omega \in \ce^k$ is closed, then
 \begin{equation*}
  e^{(n-2k)\Upsilon}Q_k^{e^{2\Upsilon}g}(\omega) = Q_k^g(\omega) +
 L_k^{n/2-k}(\Upsilon\omega) .
 \end{equation*}
\end{thm}

\noindent Here $\ce^k(w)$ denotes the space of smooth $k$-forms of conformal weight $w$ and $\ce^k := \ce^k(0)$.

The original proof of Theorem~\ref{long} uses the observation, similar
to the first construction of the GJMS operators, that each power of the ambient Hodge Laplacian commutes with
multiplication by the defining function $\sigma$ when acting on forms of the correct
homogeneity weight.
One then uses the ambient
interpretation of the tractor-$D$
operator to invariantly lift
forms and produce an invariant operator.  An alternative construction,
due to Aubry and Guillarmou~\cite{AubryGuillarmou2011}, realises the operators
$L_k^\ell$
as obstructions to regularity for a harmonic extension problem in a
Poincar\'e manifold.
Both approaches obtain the relationship of these to the $Q$-operators.

Branson and Gover pointed out many interesting interpretations of
their operators, of which we highlight one.  Let $(M,g)$ be a compact
Riemannian manifold of even dimension $n$.  For each nonnegative
integer $k \leq n/2$, the space of \defn{conformal harmonics},
\begin{equation*}
 \mathcal{H}^k := \left\{ \omega \in \ce^k \suchthat d\omega = 0 , \d Q_k\omega = 0 \right\}
\end{equation*}
is conformally invariant and finite-dimensional.
Evidently the canonical map
\begin{equation*}
 \mathcal{H}^k \ni \omega \mapsto [\omega] \in H^k(M;\mathbb{R})
\end{equation*}
onto the $k$-th de Rham group is well-defined.
It is easy to see that $\mathcal{H}^k \to H^k(M;\mathbb{R})$ is surjective in the endpoint cases $k \in \{ 0 , n/2 \}$.
Branson and Gover conjectured that this map is surjective for some class of suitably generic manifolds.

\begin{conj}
 \label{bg-conj}
 Fix an even integer $n$ and a positive integer $k < n/2$.
 For a generic
compact Riemannian $n$-manifold, the map $\mathcal{H}^k \to H^k(M;\mathbb{R})$ is surjective.
\end{conj}

Conjecture~\ref{bg-conj} represents a first step in establishing a conformally invariant Hodge theory.
If verified, it would surely have many deep applications, similar to the successes of the Hodge theorem on Riemannian manifolds.
At present, Conjecture~\ref{bg-conj} is only known to hold for
conformally Einstein manifolds~\cite{GoverSilhan2008} and for
manifolds for which $Q_k$ is nonnegative as an operator on closed
forms~\cite{AubryGuillarmou2011}.  The former result uses the
factorisation of the GJMS operators on forms. 

\subsection{Extrinsic GJMS operators} \label{extrinsic}

A classical problem of interest is to determine the conformal
invariants of a submanifold embedded in a conformal manifold $M$.  The
scattering construction~\eqref{eqn:scattering} of the GJMS operators
has been adapted, in two distinct ways, to produce natural, formally
self-adjoint, conformally invariant, linear differential operators on
these submanifolds.  These are collectively called \defn{extrinsic
  GJMS operators}.

The first construction~\cite{GoverWaldron2021}, carried out on
two-sided hypersurfaces, uses that the conformal embedding determines
(formally to a dimension-dependent order) a  metric from the conformal class on the interior of one
side of the hypersurface, namely a  singular Yamabe
metric as in \cite{AnderssonChruscielFriedrich1992}.
The second construction~\cite{CGK25},
carried out in general codimension, uses Graham and
Witten's~\cite{GrahamWitten1999} formally minimal extension of the
submanifold into a Poincar\'e manifold.  These operators are generally
different.  A notable feature of the former operators is that they
extend to all orders, regardless of the parity of the dimensions.  A
key feature of the latter operators is that they factor on minimal
submanifolds of Einstein manifolds, leading to a calculation of the
renormalised area of a minimal submanifold of a Poincar\'e--Einstein
manifold~\cite{CGKTW25}.

One can also associate \defn{extrinsic $Q$-curvatures} to the
extrinsic GJMS operators.  The integrals of these invariants on
compact submanifolds are conformal invariants that generalise the
Willmore energy.  Among other interpretations, the total extrinsic
$Q$-curvature associated to these two constructions recovers the
anomaly term $\mathcal{L}$ in the renormalised volume expansion for
singular Yamabe metrics~\cite{GoverWaldron2017,Graham2017} and the
renormalised area for the Graham--Witten minimal
submanifold~\cite{GrahamReichert2020,GrahamWitten1999}, respectively.
See~\cite{CGKTW25} for a conjectural classification of extrinsic
scalar Riemannian invariants whose integrals are conformally invariant
on compact submanifolds.

\subsection{The CR GJMS operators} \label{CR-sect}

\newcommand{\C}{\mathbb C}
\newcommand{\Co}{\mathcal C}

An important historical problem in several complex variables is to
determine when two bounded domains in $\mathbb{C}^N$, $N \geq 2$, are biholomorphic.
If the boundaries are strictly
pseudoconvex, Fefferman~\cite{Fefferman1974} showed that they are
biholomorphic if and only if their boundaries are CR equivalent.
Later, Fefferman~\cite{Fefferman1979} used this equivalence
to partially classify CR invariants.  This work directly motivated the
construction of the Fefferman--Graham ambient
metric~\cite{FeffermanGraham1985,FeffermanGraham2012}
 and an approach to the
classification of scalar conformal invariants~\cite{BaileyEastwoodGraham1994}.

The GJMS operators give an example of ideas from conformal geometry being used to further develop CR geometry.
The first such development was by Jerison and Lee~\cite{JerisonLee1987}, who found a CR analogue of the conformal Laplacian;
see Section~\ref{sec:GA/yamabe} for discussion of the CR Yamabe problem.
Gover and Graham~\cite{GoGr} constructed higher-order CR
analogues of the GJMS operators.  Indeed, in the CR setting one can
define ``bi-weights'' and find subspaces
$\mathcal{E}(w,w') \subset \mathcal{E}(w + w') \otimes \mathbb{C}$ of
the complex-valued densities when $w-w' \in \mathbb{Z}$, and
the \defn{CR GJMS operators} exist for many choices of bi-weights:

\begin{thm}
 \label{CR-GJMS}
 Fix $n \in \mathbb{N}$.
 Let $w,w' \in \mathbb{C}$ be such that $n+1+w+w'=k \in \mathbb{N}$ and $w-w' \in \mathbb{Z}$.
 If either $k \leq n+1$ or $(w,w') \not\in \mathbb{N} \times \mathbb{N}$,
 then there is a natural, formally self-adjoint, CR invariant, linear differential operator
 \begin{equation*}
  P_{w,w'}:\ce(w,w')\rightarrow \ce(w-k,w'-k)
 \end{equation*}
 on nondegenerate CR $(2n+1)$-manifolds whose principal part agrees with that of $\Delta_b^k$.
\end{thm}

Notably, the CR GJMS operators exist to all orders for generic bi-weights.
Here the ``principal part'' is determined isotropically---i.e.\ the Reeb vector field is regarded as a first-order operator.
The anisotropic principal part, which is relevant to understanding mapping properties in Folland--Stein spaces, is also known~\cite{GoGr,Ponge2008}.
These operators factor at Sasakian $\eta$-Einstein manifolds~\cite{CaseGover2020,Takeuchi2018}, the CR analogues of Einstein manifolds.

The \defn{critical} CR GJMS operators are $P:=P_{0,0}$.
In analogy with the GJMS operators on differential forms, the kernel of $P$ is generally infinite-dimensional:
it contains the space $\mathscr{P}$ of CR pluriharmonic functions.
The analogue of the $Q$-curvature operators are known as \defn{$P'$-operators}~\cite{CaseYang2012,Hirachi2013}.
These are natural linear differential operators $P' \colon \mathscr{P} \to C^\infty(M)$ such that if $\hat\theta = e^\Upsilon\theta$, then
\begin{equation*}
 e^{(n+1)\Upsilon}P_{\widehat{\theta}}'(u) = P_\theta'(u) + P_\theta(u\Upsilon) .
\end{equation*}
When $\theta$ is
pseudo-Einstein, constants are
in the kernel of $P'_\theta$, allowing one to define the \defn{$Q'$-curvature}~\cite{CaseYang2012,Hirachi2013}.  
This is a pseudohermitian scalar invariant $Q'$ such that if
$\widehat{\theta} = e^\Upsilon\theta$ is also pseudo-Einstein---or
equivalently, $\Upsilon \in \mathscr{P}$~\cite{Lee1988}---then
\begin{equation*}
 e^{(n+1)\Upsilon}Q_{\widehat{\theta}}' = Q_\theta' + P_\theta'(\Upsilon) + \frac{1}{2}P(\Upsilon^2) .
\end{equation*}
Hence $\int Q' \, \theta \wedge d\theta^n$ is independent of the choice of pseudo-Einstein contact form on a compact nondegenerate CR manifold.
This invariant is nontrivial---indeed, it recovers the Burns--Epstein invariant~\cite{BurnsEpstein1988} in CR dimension $n=1$---making it a natural generalisation of the total $Q$-curvature.
The problem of classifying all scalar pseudohermitian invariants with this invariance property remains open; 
see~\cite{Takeuchi2023} and references therein.

\section{Applications in geometric analysis}
\label{sec:ga}

The impact of the GJMS operators in geometric analysis stems from their connection to the Yamabe Problem~\cite{LeeParker1987} and its various analogues.
The \defn{$Q_{2k}$-curvature} is
\begin{equation*}
 Q_{2k} := \frac{1}{n-2k}P_{2k}(1) ;
\end{equation*}
note that $Q_2 = R/2(n-1)$ is proportional to the scalar curvature and that $Q_{n}$ is the $Q$-curvature discussed in Section~\ref{sec:alg/q}.
Finding a metric $\hg \in [g]$ of constant $Q_{2k}$-curvature is equivalent to finding a (positive if $n \not= 2k$) solution of the quasi-linear PDE
\begin{equation}
 \label{eqn:pde}
 \begin{aligned}
  P_{n}u + Q & = \lambda e^{nu} , && \text{if $n=2k$, where $\widehat g = e^{2u}g$} , \\
  P_{2k}u & = \lambda \lvert u \rvert^{\frac{4k}{n-2k}}u , && \text{if $n \not= 2k$, where $\widehat g = u^{\frac{4}{n-2k}}g$} ,
 \end{aligned}
\end{equation}
for some constant $\lambda$.
These PDEs are variational and have critical nonlinearity in the sense that their solutions are the critical points of the functionals
\begin{equation}
 \label{eqn:functional}
 \begin{aligned}
  \mathcal{F}_{2k}(u) & = \int_M \Big( uP_n u + 2Qu \Big) \dvol - \frac{2}{n}\left( \int_M Q \dvol \right) \log \int_M e^{nu} \dvol , && \text{if $2k=n$} , \\
  \mathcal{F}_{2k}(u) & = \left. \left( \int_M uP_{2k}u \dvol \right) \middle/ \left( \int_M \lvert u \rvert^{\frac{2n}{n-2k}} \dvol \right)^{\frac{n-2k}{n}} \right. , && \text{if $2k\not=n$} ,
 \end{aligned}
\end{equation}
and the continuity of these functionals on the Sobolev space $W^{k,2}(M^n)$ uses the embeddings $W^{k,2}(M^n) \hookrightarrow L^{\frac{2n}{n-2k}}(M^n)$, $n > 2k$, and $W^{n/2,2}(M^n) \hookrightarrow e^{L(M^n)}$, which are continuous but not compact.
The rest of this section requires positive signature.

\subsection{Rigidity and stability of Sharp Sobolev inequalities}
\label{sec:GA/sobolev}

A standard approach to solving~\eqref{eqn:pde} is to minimise~\eqref{eqn:functional} in $W^{k,2}(M^n)$.
As highlighted by Aubin~\cite{Aubin1976} (in the case $n > 2k$) and Chang and Yang~\cite{ChangYang1995} (in the case $n=2k$), this is facilitated by the classification of the extremals of~\eqref{eqn:functional} on the round $n$-sphere~\cite{Beckner1993}:

\begin{thm}
	\label{sharp-sobolev}
	Fix $k,n \in \mathbb{N}$ with $k \leq n/2$.
	Denote the action of the conformal group $\Conf(S^n)$ on $C^\infty(S^n)$ by
	\begin{equation*}
 		u \cdot \Phi :=
 		\begin{cases}
 			\lvert J_\Phi \rvert^{\frac{n-2k}{2n}} u \circ \Phi , & \text{if $n > 2k$} , \\
 			u \circ \Phi + \ln\lvert J_\Phi \rvert^{\frac{1}{n}} , & \text{if $n = 2k$} ,
 		\end{cases}
	\end{equation*}
	where $\lvert J_\Phi \rvert$ is the Jacobian determinant of $\Phi$.
	Then $\mathcal{F}_{2k}(u) \geq \mathcal{F}_{2k}(1)$ with equality if and only if $u = \lambda \cdot \Phi$ for some constant $\lambda \in \mathbb{R}$ and some $\Phi \in \Conf(S^n)$.
\end{thm}

The original proof of Theorem~\ref{sharp-sobolev} uses conformal
invariance and symmetric rearrangments.  Frank and Lieb~\cite{FrankLieb2012b}
gave an argument using only conformal invariance.  The benefit of
this approach is that it adapts to the CR
setting~\cite{FrankLieb2012a}, where rearrangement methods are not
available.
Case~\cite{Case2019fl} simplified their argument via the commutator
identity
\begin{equation}
	\label{eqn:commutator}
	\sum_{i=0}^n x^i[P_{2k},x^i] = k(n+2k-2)P_{2k-2}
\end{equation}
on $S^n$, where $x^i$ are the standard Cartesian coordinates.
Equation~\eqref{eqn:commutator} follows from an easy holographic computation using the first definition of the GJMS operators.

Theorem~\ref{sharp-sobolev} is a rigidity result, in that it characterises the extremal functions.
There has been a lot of work recently on obtaining a stability result;
i.e.\ a quantitative estimate on the distance of an approximate minimiser to the set of minimisers.
It is now known that if $n > 2k$, then there is a constant $C>0$ such that
\begin{equation}
	\label{eqn:bianchi-egnell}
	\mathcal{F}_{2k}(u) - \mathcal{F}_{2k}(1) \geq C \inf \left\{ \lVert \lambda u \cdot \Phi - 1 \rVert_{W^{k,2}}^2 \suchthat \lambda \in \mathbb{R}, \Phi \in \Conf(S^n) \right\} ;
\end{equation}
see the lecture notes~\cite{Frank2024} by Frank for a comprehensive discussion.

\subsection{Yamabe-type problems}
\label{sec:GA/yamabe}

The difficulties of solving~\eqref{eqn:pde} are different depending upon whether $n=2k$ or $n>2k$.
In the former case, the key difficulty is that $\mathcal{F}_n$ need not be bounded below, so critical points cannot always be constructed via minimisation methods.
In the latter case, the key difficulty is that one cannot apply elliptic regularity to conclude smoothness unless $u>0$, which is difficult given the lack of a maximum principle for the GJMS operators.
We discuss these cases separately.

Suppose that $n=2k$.
When $k=1$, a solution of~\eqref{eqn:pde} gives a proof of the Uniformisation Theorem on compact surfaces.
A PDE proof was given by Osgood, Phillips, and Sarnak~\cite{OsgoodPhillipsSarnak1988}.
Their result was later generalised to dimension four by Chang and Yang~\cite{ChangYang1995} and to general dimensions by Brendle~\cite{Brendle2003}:

\begin{thm}
	\label{critical-yamabe}
	Let $(M^n,g)$ be a compact even-dimensional Riemannian manifold such that $P_n \geq 0$, $\ker P_n = \mathbb{R}$, and $\int_M Q \dvol < \int_{S^n}Q \dvol$.
	There is there is smooth minimiser $u$ of $\mathcal{F}_n$.
	In particular, $e^{2u}g$ has constant $Q$-curvature.
\end{thm}

The assumptions of Theorem~\ref{critical-yamabe} ensure that
$\mathcal{F}_n$ can be minimised.  Djadli and
Malchiodi~\cite{DjadliMalchiodi2008} pioneered the use of minimax
methods to find critical points of $\mathcal{F}_n$ under weaker
hypotheses; see~\cite{LiLiLiu2012,Ndiaye2007} for recent developments.
Note, however, that Equation~\eqref{Qtrans} implies that $\int Qu
\dvol$ is conformally invariant for all $u \in \ker P_n$, leading to
constraints on $Q$ when $\ker P_n\neq \ker d $~\cite{Gover2010}.

Suppose now that $n > 2k$.
When $k=1$, a solution of~\eqref{eqn:pde} solves the Yamabe Problem;
see the survey by Lee and Parker~\cite{LeeParker1987}.
When $k=2$, one of the most general existence results is due to Gursky, Hang, and Lin~\cite{GurskyHangLin2016}.
In the general case, Mazumdar~\cite{Mazumdar2016} found two global properties that are sufficient to solve~\eqref{eqn:pde} by minimising the \defn{$Q_{2k}$-Yamabe constant}
\begin{equation*}
	Y_{2k}(M^n,[g]) := \inf \left\{ \mathcal{F}_{2k}(u) \suchthat 0 \not= u \in C^\infty(M) \right\} :
\end{equation*}

\begin{thm}
	\label{subcritical-yamabe}
	Let $(M^n,g)$ be a compact Riemannian manifold and let $k < n/2$ be a positive integer.
	Suppose that
	\begin{enumerate}
		\item [(A)] $0 < Y_{2k}(M^n,[g]) < Y_{2k}(S^n,[d\theta^2])$, and
		\item [(M)] for every $u \in C^\infty(M)$, if $P_{2k}u \geq 0$, then $u > 0$ or $u = 0$.
	\end{enumerate}
	Then there is a smooth positive minimiser of $Y_{2k}(M,[g])$.
\end{thm}

\noindent Here $Y_{2k}(S^n,[d\theta^2])$ is the $Q_{2k}$-Yamabe constant of the round $n$-sphere.
Its value is computed by Theorem~\ref{sharp-sobolev}.

The upper bound in Property~(A) ensures that minimising sequences for $Y_{2k}(M,[g])$ converge in $W^{k,2}$, generalising an observation of Aubin~\cite{Aubin1976}.
Property~(M) and the lower bound in Property~(A) together ensure that minimisers are positive.

For general $k$, Property~(A) was verified by Qing and Raske~\cite{QingRaske2006} for locally conformally flat manifolds with Poincar\'e exponent $\delta < \frac{n-2k}{2}$, and by Mazumdar and V\'etois~\cite{MazumdarVetois2024} for manifolds that are not locally conformally flat provided $n \geq 2k+4$;
the latter argument uses Juhl's formula~\eqref{eqn:juhl-formula}.
Hang and Yang~\cite{HangYang2016} verified Property~(M) when $k=2$ under the assumption that $Y_2(M,[g])>0$ and there is a metric $\hg \in [g]$ with $Q_4^{\hg} \gneqq 0$.
Inspired by this work, Andrade, Piccione, and Wei~\cite{AndradePiccioneWei2023} conjectured a sufficient condition for $P_{6}$ to satisfy Property~(M).
We formulate their conjecture for general orders:

\begin{conj}
	\label{smp}
	Let $(M^n,g)$ be a compact Riemannian manifold and let $k < n/2$ be a positive integer.
	If $Y_{2j}(M,[g])>0$ for each $j \in \{ 1, \dotsc, k-1 \}$ and if there is a metric $\hg \in [g]$ with $Q_{2k}^{\hg} \gneqq 0$, then $P_{2k}$ satisfies Property~(M).
\end{conj}

There has also been work on the (lack of) uniqueness of solutions to~\eqref{eqn:pde}.
Two trivial sources of the lack of compactness
are the action of the conformal group and multiplication by a constant.
All statements below are made modulo these two actions.

Uniqueness for conformally Einstein metrics was proved by Obata~\cite{Obata1971} when $k=1$ and by V\'etois~\cite{Vetois2024} when $k=2$ and $Y_2(M,[g])\geq0$.
It is unknown if a similar conclusion is true for larger $k$.

Schoen~\cite{Schoen1989} exhibited nonuniqueness when $k=1$ by
studying product metrics on $S^1 \times S^{n-1}$.  Bettiol and
Piccione~\cite{BettiolPiccione2018} generalised this to manifolds with
$Y_2(M,[g])>0$ whose fundamental group has infinite profinite
completion.  Their conclusion remains true for all $k < n/2$ when
Theorem~\ref{subcritical-yamabe} can be
applied~\cite{AndradeCasePiccioneWei2023}.  Moreover, Brendle and
Marques~\cite{Brendle2008,BrendleMarques2009} proved the remarkable fact that there are nontrivial
compact manifolds of dimension $n \geq 25$ for which the set of
solutions to~\eqref{eqn:pde} with $k=1$ is noncompact; this dimension
bound is sharp~\cite{KhuriMarquesSchoen2009}.  When $k=2$, dimension
$n=25$ is \emph{also} the minimal dimension for which the solution set
can be noncompact~\cite{GongKimWei2025,WeiZhao2013}.  The case $k \geq
3$ is open.

The quantitative stability of the $Q_{2k}$-Yamabe constant has also been studied.
The analogue of~\eqref{eqn:bianchi-egnell} holds for $k=1$, though one must replace the exponent $2$ on the right-hand side by $2+\gamma$~\cite{EngelsteinNeumayerSpolaor2022}.
Moreover, this is sharp with $\gamma \in [0,2]$, as can seen from the bifurcation points in Schoen's analysis of $S^1 \times S^{n-1}$~\cite{Frank2022}.
Quantitative stability is also known for $k \geq 2$ when the hypotheses of Theorem~\ref{subcritical-yamabe} are satisfied~\cite{AndradeKonigRatzkinWei2024}.

\subsection{Fractional Yamabe problems}
\label{sec:GA/nonlocal}

There are natural, formally self-adjoint, conformally covariant,
linear pseudodifferential operators $P_{2\gamma}$ with leading-order
term $\Delta^\gamma$; see Section~\ref{factor}.  Gonz\'alez and
Qing~\cite{GonzalezQing2013} initiated the study of the corresponding
Yamabe-type problem; see~\cite{KimMussoWei2021,MayerNdiaye2024} and
references therein for the most recent progress.

When $\gamma=1/2$, the fractional Yamabe problem coincides~\cite{GonzalezQing2013} with the a boundary Yamabe problem introduced by Escobar~\cite{Escobar1992}.
More generally, \emph{all} of the fractional Yamabe problems correspond to (weighted) boundary Yamabe problems involving higher-order GJMS operators~\cite{CaseChang2016}.
However, in many cases these boundary Yamabe problems are currently only defined on compactifications of Poincar\'e--Einstein manifolds;
see~\cite{Case2018,FlynnLuYang2023} and references therein.

\subsection{CR Yamabe-type problems}
\label{sec:GA/cr-yamabe}

Jerison and Lee~\cite{JerisonLee1987,JerisonLee1988,JerisonLee1989} posed and studied the Yamabe-type problem for the second-order CR GJMS operator:
They constructed minimisers of the CR analogue of the functional~\eqref{eqn:functional} under the CR analogue of the assumption $Y_2(M^n,\cc) < Y_2(S^n,[d\theta^2])$, verified this assumption when $n \geq 2$ and $(M^{2n+1},T^{1,0})$ is not locally CR equivalent to the standard CR sphere, and classified the extremals on the standard CR sphere.
As previously noted, Frank and Lieb~\cite{FrankLieb2012a} gave an alternative proof of this classification via factorisation which is valid for all CR GJMS operators; see also \cite{Yan2023}.
The resolution of the CR Yamabe Problem when $n=1$, which requires embeddability and a positive mass theorem, was accomplished by Cheng, Malchiodi, and Yang~\cite{ChengMalchiodiYang2017};
see~\cite{ChengMalchiodiYang2023} for counterexamples without these assumptions.
Yan~\cite{Yan2023} proved the CR analogue of Theorem~\ref{subcritical-yamabe} when $k=2$, but the higher-order CR Yamabe Problem remains otherwise open.

The CR analogue of Theorem~\ref{critical-yamabe}, which prescribes the $Q'$-curvature to be constant, was accomplished by Case, Hsiao, and Yang~\cite{CaseHsiaoYang2019} when $n=1$.

\subsection{Fully nonlinear generalisations}
\label{sec:GA/fully-nonlinear}

The Yamabe problem has also inspired work on fully nonlinear problems, especially those stated in terms of the $\sigma_k$-curvatures;
i.e.\ the $k$-th elementary symmetric functions of the Schouten tensor
\begin{equation*}
 P := \frac{1}{n-2}\left( \Ric - \frac{R}{2(n-1)}g \right) .
\end{equation*}
The problem of conformally prescribing $\sigma_k$ to be constant is variational when $k=2$ or $g$ is locally conformally flat~\cite{BransonGover2008}, and in this setting the problem is solved~\cite{GeWang2006,GuanWang2004,ShengTrudingerWang2007}.

The holographically-defined renormalised volume coefficients are variational scalar invariants which equal the $\sigma_k$-curvatures when the latter are variational;
in general, their difference involves lower-order terms~\cite{Graham2009}.
Their conformal transformation law can be written in terms of a conformally invariant multilinear differential operator~\cite{CaseLinYuan2022}.
In the special case of the $\sigma_2$-curvature, this operator has a simple holographic formula that was used to prove the rigidity~\cite{Case2019fl} and stability~\cite{FrankPeteranderl2025} of the analogue of Theorem~\ref{sharp-sobolev}.

\section*{Acknowledgements}
We thank Mike Eastwood for helpful comments on the history of the GJMS operators.

JSC was partially supported by a Simons Foundation Collaboration Grant for Mathematicians and by the National Science Foundation under Award No.\ DMS-2505606.
RG gratefully acknowledges support from the Royal
Society of New Zealand via Marsden Grants 19-UOA-008 and 24-UOA-005.

\bibliographystyle{abbrv}
\bibliography{bib}
  
\end{document}